\documentclass{amsart}
\usepackage{amssymb}
\usepackage{amscd}

\newtheorem{Theorem}{Theorem}[section]
\newtheorem{Lemma}[Theorem]{Lemma}
\newtheorem{Corollary}[Theorem]{Corollary}
\newtheorem{Proposition}[Theorem]{Proposition}
\theoremstyle{definition}

\newtheorem{Example}[Theorem]{Example}
\theoremstyle{remark}
\newtheorem{Remark}{Remark}

\font\sy=cmsy10
\font\ym=msbm10
\newcommand{\cC}{{\hbox{\sy C}}}
\newcommand{\cD}{{\hbox{\sy D}}}
\newcommand{\cE}{{\hbox{\sy E}}}
\newcommand{\cH}{{\hbox{\sy H}}}
\newcommand{\cK}{{\hbox{\sy K}}}
\newcommand{\cM}{{\hbox{\sy M}}}
\newcommand{\cT}{{\hbox{\sy T}}}
\newcommand{\cV}{{\hbox{\sy V}}}
\newcommand{\C}{{\text{\ym C}}}

\newcommand{\R}{\text{\ym R}}

\newcommand{\Z}{\text{\ym Z}}
\newcommand{\trace}{\hbox{\rm tr}}
\newcommand{\Hom}{\hbox{\rm Hom}}
\newcommand{\End}{\hbox{\rm End}}
\title[]
{Fiber Functors on Temperley-Lieb Categories}
\author[Yamagami Shigeru]{Shigeru Yamagami}

\begin{document}
\maketitle   
\begin{center}
Department of Mathematical Sciences
\end{center}
\begin{center}
Ibaraki University 
\end{center}
\begin{center} 
Mito, 310-8512, JAPAN 
\end{center}    

\begin{abstract}
Fiber functors on Temperley-Lieb categories are determined 
with the help of classification results on non-degenerate 
bilinear forms. The case of unitary fiber functors is also 
investigated.
\end{abstract}
\bigskip
%\baselineskip=18pt
% \rotatebox{180}{C}
% \rotatebox{60}{Q}
\section{Introduction}
Tensor categories are often realized 
as categories of finite-dimensional comodules of Hopf algebras, 
which can be seen 
in other ways as giving a faithful exact monoidal functor 
from an abelian tensor category 
into that of finite-dimensional vector spaces. 
The converse implication holds if one imposes rigidity 
on the tensor category in question; 
given a faithful exact monoidal functor $F$ 
from a rigid abelian tensor category $\cT$ 
into the tensor category of finite-dimensional vector spaces, 
there exists a Hopf algebra $A$ such that $\cT$ is monoidally 
equivalent to the tensor category of 
finite-dimensional $A$-comodules with the functor identified with 
the one forgetting $A$-coactions (\cite{S,U}). 

An exact monoidal functor from an abelian tensor category into 
the tensor category $\cV ec$ of 
finite-dimensional vector spaces is usually called 
a fiber functor. 

Tensor categories we shall work with are the so-called 
Temperley-Lieb categories 
$\cK_d$ ($d \in \C^\times = \C \setminus \{ 0\}$): 
$\cK_d$ is a rigid 
tensor category associated with planar tangles, whose 
objects are parametrized natural numbers $n = 0, 1, 2, \dots$ and 
$\End(n) = \Hom(n,n)$ is in an obvious way identified with 
the $n$-th Temperley-Lieb algebra (generated by unit and elements 
$h_1,\dots, h_{n-1}$ satisfying $h_j^2 = d h_j$, 
$h_ih_j = h_jh_i$ ($|i-j| \geq 2)$ and $h_ih_jh_i = h_i$ ($|i-j| = 1$)).
The Temperley-Lieb category $\cK_d$ admits a compatible C*-structure 
if and only if $|d| \geq 2$ with $d \in \R$ or 
$d = \pm 2\cos (\pi/l)$ with $l = 3, 4, \dots$. In that case, 
the compatible C*-structure is unique up to C*-monoidal equivalences 
(see, for example, \cite{CDA} and references therein). 

Although the Temperley-Lieb category $\cK_d$ is never abelian, 
the idempotent-completion of $\cK_d$ after adding finite-direct sums 
(first add direct sums and then add subobjects) 
turns out to be semisimple and therefore abelian unless  
$d = q + q^{-1}$ with $q^2$ a non-trivial root of unity.
Moreover, any monoidal functor $F: \cK_d \to \cV ec$ is uniquely 
extended to this abelianized Temperley-Lieb category, 
which is automatically exact by semisimplicity. 

With these backgrounds in mind, we shall refer to a monoidal 
functor $\cK_d \to \cV ec$ simply 
as a \textbf{fiber functor} in what follows. 
When $\cK_d$ is a C*-tensor category, we can naturally talk 
about the unitarity of fiber functors; a monoidal functor 
$F: \cK_d \to \cH ilb$ ($\cH ilb$ denotes the C*-tensor category 
of finite-dimensional Hilbert spaces) is called 
a \textbf{unitary fiber functor} 
if $F$ preserves the *-operation. 
For a unitary fiber functor, the reconstructed Hopf algebra 
is naturally regarded as defining a compact quantum group $G$ 
with its representation category (the Tannaka dual of $G$) 
monoidally 
equivalent to the C*-tensor category $\cK_d$. 

Our main concern here is how fiber functors are determined 
on Temperley-Lieb categories and we shall present 
a classification result  
with the help of similar works on 
non-degenerate bilinear forms 
(see \cite{Wa} and references therein). 
By Saavedra Rivano-Ulbrich's theorem on Tannaka-Krein duality, 
this produces 
mutually non-isomorphic Hopf algebras 
with the isomorphic representation 
categories as far as $d$ is generic 
in the sense that $\cK_d$ is semisimple.
Since the relevant Hopf algebras turn out to be the ones 
introduced by M.~Dubois-Violette and G.~Launer (\cite{DVL}), 
we particularly obtain 
an isomorphism classification of Hopf algebras 
in that class. 

When $|d| \geq 2$ with $d \in \R$ and fiber functors are restricted 
to be unitary, we also give their complete descriptions, 
resulting in 
a classification of compact quantum groups 
whose Tannaka duals are isomorphic to 
$\cK_d$ as monoidal categories. 
In this case, the relevant Hopf algebras are identified with   
the free orthogonal quantum groups 
investigated by T.~Banica and S.Z.~Wang
(\cite{Ba}, \cite{DW}, \cite{Wan}),  
and the present analysis enables us to get access to 
some of their results from the viewpoint of 
tensor categories.

\section{Fiber Functors}
Let $F: \cK_d \to \cV ec$ be a fiber functor. 
By definition, it consists of a linear functor $F$ together with 
a natural family of isomorphisms 
$\mu_{m,n}: F(X^m)\otimes F(X^n) \to F(X^{m+n})$ which makes 
the following diagram commutative 
\[
\begin{CD}
F(X^l)\otimes F(X^m)\otimes F(X^n) @>{1\otimes \mu_{m,n}}>> 
F(X^l)\otimes F(X^{m+n})\\
@V{\mu_{l,m}\otimes 1}VV @VV{\mu_{l,m+n}}V\\
F(X^{l+m})\otimes F(X^n) @>>{\mu_{l+m,n}}> F(X^{l+m+n})
\end{CD}.
\]
%We can then show that there exist well-defined isomorphisms 
Then isomorphisms $F(X^n) \to F(X)^{\otimes n}$ obtained 
as repetitions of $\mu$'s are identical and define 
a single isomorphism (the coherence theorem), which satisfies 
the commutativity of the diagram 
\[
\begin{CD}
F(X^m)\otimes F(X^n) @>{\mu_{m,n}}>> F(X^{m+n})\\
@VVV @VVV\\
F(X)^{\otimes m}\otimes F(X)^{\otimes n} 
@= F(X)^{\otimes (m+n)}
\end{CD}.
\]
In other words, the functor $F$ is monoidally equivalent 
to strict one. 
Then we see that the functor $F$ in its strict version is 
specified by the choice of 
% By a kind of coherence thereom, we may assume that $F$ is strict 
% in the sense that $F(X^m)\otimes F(X^n) = F(X^{m+n})$. 
% In particular, $F(X^n) = F(X)^{\otimes n}$ and $F$ is uniquely 
% determined by the choice of 
a bilinear form $F(\epsilon): F(X)\otimes F(X) \to \C = F(I)$, 
which is non-degenerate by rigidity 
(cf.~\cite[Lemma 6.1]{CDA}, tensor categories being 
assumed strict without qualifications there).

Conversely, given a non-degenerate bilinear form 
$\cE: V\otimes V \to \C$ on a finite-dimensional vector space $V$, 
we can determine the copairing $\cD: \C \to V\otimes V$ by 
the rigidity identity and, if the dimension relation 
$\cE \cD = d$ is satisfied furthermore, 
we can recover the strict monoidal functor 
$F: \cK_d \to \cV ec$ so that $F(X) = V$ and $F(\epsilon) = \cE$. 

Two fibre functors $F$, $G: \cK_d \to \cV ec$ are said to be 
(monoidally) equivalent if there is a natural equivalence 
$\{ t_n: F(X^n) \to G(X^n) \}$ fulfilling the commutativity of 
the diagram 
\[
\begin{CD}
F(X^m)\otimes F(X^n) @>{\mu_F}>> F(X^{m+n})\\
@V{t_m\otimes t_n}VV @VV{t_{m+n}}V\\
G(X^m)\otimes G(X^n) @>>{\mu_G}> G(X^{m+n})
\end{CD}.
\]

By using the natural isomorphisms 
$F(X^n) \to F(X)^{\otimes n}$ and 
$G(X^n) \to G(X)^{\otimes n}$ explained above, 
the natural equivalence $\{ t_n\}$ is uniquely determined 
by the initial isomorphism $t_1: F(X) \to G(X)$ which makes 
the following diagram commutative. 
\[
\begin{CD}
F(X)\otimes F(X) @>{t_1\otimes t_1}>> G(X)\otimes G(X)\\
@V{F(\epsilon)}VV @VV{G(\epsilon)}V\\
\C @= \C
\end{CD}
\]
where $F(I) = \C = G(I)$. 

Conversely, starting with a linear isomorphism 
$t_1: F(X) \to G(X)$ making the above diagram commutative, 
the choice $t_n = (t_1)^{\otimes n}$ gives 
a monoidal equivalence between $F$ and $G$.
 
Thus, if we take the canonical column vector spaces $\C^N$ 
with the canonical basis $\{ e_i\}$ 
as representatives of finite-dimensional vector spaces and 
describe a non-degenerate bilinear form $\cE$ on $\C^N$ by 
the associated invertible matrix $E \in GL_N(\C)$ 
\[
\cE(e_i\otimes e_j) = E_{i,j},
\]
then the accompanied copairing 
$\cD = \cD(1) \in \C^N\otimes \C^N$ is given by 
\[
\cD = \sum_{i,j} D_{i,j} e_i\otimes e_j,
\]
with the matrix $D = \{ D_{i,j}\}$ characterized 
as the inverse of $E$. 
The dimension condition $d = \cE \cD$ then takes the form 
%(the associated dimension is given by the formula)
\[
d = \trace({}^tE^{-1} E) = \trace({}^tE E^{-1}).
\]

%% Thus natural equivalence classes of fiber functors on $\cK_d$ are in one-to-one 
%% correspondance 
%% with orbits $O(G)$ in the set $GL_N(\C)$ under the (right) action 
%% of the group $GL_N(\C)$ defined by 
%% \[
%% G \mapsto {}^tT GT.
%% \]

% If the function $d$ is defined on $GL_N(\C)$ by the left hand side 
% of the above formula, it is an invariant of the $GL_N(\C)$-action 
% on itself defined by $a.g = ag\,{}^ta$, i.e., 
% $d(ag\,{}^ta) = d(g)$. 

Let $M_N(d)$ be the set of such matrices in $GL_N(\C)$.  
% Two matrices $G$, $G' \in M_N(d)$ then gives rise to 
% naturally equivalent fibre functors if and only if $O(G) = O(G')$
Since ${}^tT E T$ belongs to $M_N(d)$ 
whenever $E \in M_N(d)$ and $T \in GL_N(\C)$, 
fiber functors are completely classified up to monoidal equivalences 
by the orbit space 
$M_N(d)/GL_N(\C)$, where the right action of $GL_N(\C)$ on 
the set $M_N(d)$ is defined by 
\[
E \mapsto {}^tT E T, 
\qquad 
T \in GL_N(\C),\ E \in M_N(d).
\]

This kind of classification problem 
on transposed similarity has been investigated in 
various contexts and we know 
a complete structural analysis (\cite{Wa, R, Ha}). 
To describe the relevant results, 
we introduce some notations. 

Let $\Theta: GL_N(\C) \to GL_N(\C)$ be the map defined by 
$\Theta(E) = (E^{-1}) ({}^tE)$ for $E \in GL_N(\C)$. 
Then the relation $\Theta({}^tTET) = T^{-1}\Theta(E)T$ 
reveals that $\Theta$ is 
equivariant with respect to 
the transposed similarity action on the domain and 
the ordinary similarity action on the range. 
The similarity class of an matrix $M$ is 
in turn completely described by the multiplicity sequence 
$\mu_M(z) = (\mu_M^{(1)}(z), \mu_M^{(2)}(z), \dots)$ 
with $z \in \C$, 
where $\mu_M^{(k)}(z)$ denotes 
the multiplicity of $z$-Jordan block of size $k$ in $M$. 
Thus $\mu_M(z) = 0$ unless $z$ is an eigenvalue of $M$ and 
the sequence $\mu_M(z)$ 
vanishes after the $N$-th component. 

\begin{Theorem}[\cite{Wa}, cf.~also \cite{R, Ha}]~ 
  \begin{enumerate}
  \item 
Two invertible matrices $E$ and $E'$ are equivalent 
relative to the transposed similarity 
if and only if $\Theta(E)$ and $\Theta(E')$ are conjugate 
in $GL_N(\C)$, i.e., 
we can find $T \in GL_N(\C)$ such that 
$\Theta(E') = T^{-1} \Theta(E) T$. 
\item 
An invertible matrix $M \in GL_N(\C)$ belongs to 
the image of $\Theta$ if and only if 
(i) $\mu_M(z) = \mu_M(z^{-1})$ for $z \in \C^\times$, 
(ii) the multiplicity 
$\mu_M^{(k)}(1)$ is even for even $k$ and (iii) the multiplicity 
$\mu_M^{(k)}(-1)$ is even for odd $k$.
  \end{enumerate}
\end{Theorem}

For a multiplicity sequence $\mu = (\mu^{(1)}, \mu^{(2)}, \dots)$, 
define its total multipicity $|\mu|$ by 
\[
|\mu| = \sum_{k \geq 1} k \mu^{(k)}.
\]
Note that 
the function $|\mu_M(z)|$ is in one-to-one correspondance with 
the characteristic polynomial 
\[
\det(\lambda I - M) = \prod_z (\lambda - z)^{|\mu_M(z)|}. 
\]

In particular, the trace $\trace(\Theta(E))$ of 
$\Theta(E)$ is an invariant of 
the transposed similarity of $E$. 

If we combine this result with our previous arguments, the following 
is obtained. 

\begin{Theorem}
There is a one-to-one correspondance between 
isomorphism classes of fiber functors 
on the Temperley-Lieb category $\cK_d$ and multiplicity functions 
$\mu(z) = (\mu^{(1)}(z), \mu^{(2)}(z), \dots)$ 
for $z \in \C^\times$ satisfying 
the conditions (i)-(iii) in the above theorem and 
\[
d = \sum_{z \in \C^\times} |\mu(z)|\, z.
\]
\end{Theorem}

If the multiplicity function $\mu(z)$ meets the condition that 
$\mu^{(k)}(\pm 1)$ is even except for $\mu^{(1)}(1)$, 
then we can describe 
a matrix $E$ in the following simple way; 
let $\{ q_i \}$ be a half of the spectral set 
$\sigma = \{ z \in \C^\times; \mu(z) \not= 0\}$ 
satisfying 
$\sigma \setminus \{ \pm 1\} = \{ q_i, q_i^{-1} \}$. 
Choose an invertible (square) matrix $Q$ so that 
\[
m_Q(z) = 
\begin{cases}
\mu(q_i) &\text{if $z = q_i$ for some $i$.}\\
\mu(-1)/2 &\text{if $z = -1$,}\\
0 &\text{otherwise}
\end{cases}
\]
and set 
\[
E = 
\begin{pmatrix}
0 & I & 0\\
Q^{-1} & 0 & 0\\
0 & 0 & I 
\end{pmatrix}.
\]
Then it is immediate to see that the multiplicity function of 
\[
\Theta(E) = 
\begin{pmatrix}
Q & 0 & 0\\
0 & {}^tQ^{-1} & 0\\
0 & 0 & I
\end{pmatrix}
\]
is given by the initial function $\mu(z)$.

\begin{Example}
The generic orbits in $GL_2(\C)$ are represented by 
\[
E_q = 
\begin{pmatrix}
0 & 1\\
-q^{-1} & 0
\end{pmatrix}
\qquad
\text{with}
\quad
\Theta(E_q) = - 
\begin{pmatrix}
q & 0\\
0 & q^{-1}
\end{pmatrix}
\]
for $q \in \C^\times$.

The stabilizer at $E_q$ is given by 
\[
\left\{ 
  \begin{pmatrix}
    t & 0\\
    0 & t^{-1}
  \end{pmatrix}; t \in \C^\times 
\right\} 
\cong \C^\times 
\]
for $q \not= \pm 1$, 
\[
\left\{ 
  \begin{pmatrix}
    s & 0\\
    0 & s^{-1}
  \end{pmatrix},
  \begin{pmatrix}
    0 & t\\
t^{-1} & 0
  \end{pmatrix}
; s, t \in \C^\times 
\right\}
\cong \C^\times \rtimes \Z_2
\]
for $q=-1$ and $SL_2(\C)$ for $q=1$. 

There remains one more orbit, which is represented by 
\[
E = 
\begin{pmatrix}
1 & 1\\
-1 & 0
\end{pmatrix}
\qquad
\text{with}
\quad
\Theta(E) = 
\begin{pmatrix}
-1 & 0\\
2 & -1  
\end{pmatrix}, 
\]
where the stabilizer is given by 
\[
\left\{
\begin{pmatrix}
\pm 1 & z\\
0 & \pm 1
\end{pmatrix}; 
z \in \C
\right\}
\cong \C\rtimes \Z_2.
\]

Thus there are unique fiber functors when $d \not= -2$ 
while we have two non-isomorphic ones for 
the quantum dimension $d = -2$.
\end{Example}

By utilizing these explicit descriptions, we can write down the 
associated Hopf algebra $A$ as well. %in a concrete manne. 
If we rename the generator $\{ a_{ij} \}_{1 \leq i,j \leq 2}$ 
of the associated Hopf algebra $A$ by 
\[
\begin{pmatrix}
a_{11} & a_{12}\\
a_{21} & a_{22}
\end{pmatrix}
= 
\begin{pmatrix}
a & b\\
c & d
\end{pmatrix}
\]
and apply the description in the appendix, 
then the quantum group $\text{SL}_q(2,\C)$ is recovered 
(see \cite{Kas} for more information on $\text{SL}_q(2,\C)$); 
$A$ is the unital algebra $A_q$ 
generated by $a$, $b$, $c$ and $d$ with 
the relations 
\begin{gather*}
ab = q^{-1}ba, ac = q^{-1}ca, bd = q^{-1}db, 
cd = q^{-1}dc,\\
 bc = cb, 
ad - q^{-1}bc = da - qbc = 1
\end{gather*}
and the coproduct is given by 
\[
\Delta
\begin{pmatrix}
a & b\\
c & d
\end{pmatrix}
= 
\begin{pmatrix}
a & b\\
c & d
\end{pmatrix}
\otimes 
\begin{pmatrix}
a & b\\
c & d
\end{pmatrix}. 
\]

In view of these commutation relations, 
we notice that the (polynomial) function algebra on 
the matrix group $SL(2,\C)$ comes out 
from the orbit specified by $q=1$. 

Next consider the singular orbit, which is represented by 
matrices
\[
\begin{pmatrix}
\tau & 1\\
-1 & 0
\end{pmatrix}, 
\qquad 0 \not= \tau \in \C.
\]
Although they belong to the same orbit, this redundancy 
will be helpful in identifying the associated Hopf algebra $A$. 
Keeping the notation for generators of $A$, 
the Hopf algebra $A$ turns out to 
be the quantum group dealt with in \cite{W3}; 
$A$ is the unital algebra generated by $a$, $b$, $c$ and $d$ with 
the relations 
\begin{gather*}
\tau a^2 - ba + ab = \tau 1, 
\quad
\tau d^2 - bd + db = \tau 1,\\
\tau ac - bc + ad = 1, 
\quad
\tau dc - bc + da = 1,\\
\tau ca - da + cb = -1,
\quad
\tau cd - ad + cb = -1,\\ 
\tau c^2 - dc + cd = 0,
\quad
\tau c^2 - ac + ca = 0
\end{gather*}
with the coproduct given by the same formula as above. 
In particular, we see that all these Hopf algebras are 
isomorphic irrelevant to the choice of $0 \not= \tau \in \C$. 

\section{Unitary Fiber Functors}

Let $\pm d \geq 2$ and furnish the Temperley-Lieb category 
$\cK_d$ with a canonical C*-structure so that it is a C*-tensor 
category (see \cite[Proposition~6.2]{CDA}). 

In general, given a C*-tensor category $\cT$, 
a monoidal functor $F$ from $\cT$ into the C*-tensor category 
$\cH ilb$ of finite-dimensional Hilbert spaces 
is called a \textbf{unitary fiber functor} if the multiplicativity 
isomorphisms $m_{X,Y}: F(X)\otimes F(Y) \to F(X\otimes Y)$ 
are unitary and the linear maps
$F: \Hom(X,Y) \to \Hom(F(X),F(Y))$ on hom-spaces are *-preserving.

Two unitary fiber functors $F$, $G: \cT \to \cH ilb$ are said to be 
unitarily equivalent if there is a monoidal natural equivalence 
$\{ \varphi_X: F(X) \to G(X) \}$ given by unitary maps. 

For the C*-tensor category $\cK_d$, giving a unitary fiber functor 
is equivalent to specifing a finite-dimensional Hilbert space $V$ 
and linear maps 
\[
F(\epsilon): V\otimes V \to \C
\quad
\text{and}
\quad 
F(\delta): \C \to V\otimes V
\]
satisfying the rigidity relations and the unitarity condition 
$F(\epsilon)^* = (d/|d|) F(\delta)$. 

Recall here some categorical operations on Hilbert spaces. Given 
a Hilbert space $V$, denote by $\overline V$ the conjugate Hilbert 
space of $V$. By the `self-duality' of Hilbert spaces, 
$\overline V$ is identified with the dual Hilbert space of $V$. 
As a notation, 
we use $\overline v$ to stand for the associated vector 
in $\overline V$ so that $\overline v$ defines 
a linear functional by 
$\langle \overline v,w\rangle = (v|w)$ for $v, w \in V$ 
(the inner product $(v|w)$ being linear in $w$ by our convention). 
Given a (bounded) linear map $f: V \to W$ between Hilbert spaces, 
we denote by 
$\overline f: \overline V \to \overline W$ the linear map 
defined by ${\overline f}{\overline v} = \overline{f(v)}$, 
which is referred to as the conjugation of $f$.  
The operation of conjugation commutes with that of 
taking hermitian conjugate. 
The transposed operation is then defined by 
${}^tf = (\overline f)^* = \overline{f^*}$, which is a linear map 
$\overline W \to \overline V$.
Thus, the three operations ${}^tf$, $f^*$ and $\overline f$ are 
mutually commutative. Moreover, for invertible linear maps, 
these operations preserve inverses. 

Now, given a non-degenerate bilinear form 
$\cE: V\otimes V \to \C$, define a linear map 
$\Phi: V \to \overline V$ by 
\[
\cE(v\otimes v') = (\overline v| \Phi v').
\]

Let $\cD: \C \to V\otimes V$ be the associated coparing. 
If we identify $\cD$ with a vector $\cD(1)$ in $V\otimes V$, 
then an expression
\[
\cD = \sum_j v_j\otimes \delta_j
\]
with $\{ v_j\}$ an orthonormal basis in $V$ satisfies 
the rigidity relation 
if and only if $\delta_j = {}^t\Phi^{-1} \overline{v_j}$, i.e., 
\[
\cD = \sum_j v_j \otimes {}^t\Phi^{-1} \overline{v_j}.
\]
As a benefit of this expression, we have 
\[
d = \cE \circ \cD = \sum_j 
(\overline{v_j}|\Phi\, {}^t\Phi^{-1} \overline{v_j}) 
= \trace(\Phi\, {}^t\Phi^{-1}).
\]

The unitarity condition $\cE^* = (d/|d|) \cD$ is then 
equivalent to 
\begin{align*}
(\overline{v}|\Phi w) &= \cE(v\otimes w) 
= (\cE^*| v\otimes w)
= \frac{d}{|d|} (\cD| v\otimes w)\\
&= \frac{d}{|d|} 
\sum_j (v_j\otimes {}^t\Phi^{-1} \overline{v_j} | v\otimes w)\\
&= \frac{d}{|d|} \sum_j (\overline{v}|\overline{v_j}) 
(\overline{v_j} | \left( {}^t\Phi^{-1} \right)^* w)\\
&= \frac{d}{|d|} (\overline v | {\overline \Phi}^{-1} w), 
\end{align*}
i.e., the invertible map $\Phi: V \to \overline V$ should satisfy 
the relation  
\[
\Phi^{-1} = \frac{d}{|d|} \overline \Phi,
\]
which is equivalent to ${}^t\Phi^{-1} = (d/|d|) \Phi^*$. 
Note here that this implies the relation  
\[
\trace(\Phi \Phi^*) 
= \frac{|d|}{d} \trace(\Phi\,{}^t\Phi^{-1}) 
= |d|.
\]

% Choosing an orthonormal basis $\{ v_i\}$ in $V$, the bilinear form 
% $F(\epsilon)$ is described by the matrix 
% $(F_{ij} = F(\epsilon)(v_i\otimes v_j))$ with $F(\delta)$ determined by 
% the hook identities as 
% \[
% F(\delta)(1) = \sum (F^{-1})_{ij} v_i\otimes v_j
% \]
% and the unitarity condition takes the form 
% \[
% F^{-1} = \frac{d}{|d|} \overline F.
% \]
% Thus, the dimension formula 
% \[
% d = \trace({}^tF^{-1}F) = \frac{d}{|d|} \trace(F^*F)
% \]
% forces the relation $\trace(F^*F) = |d|$. 

Conversely, starting with an invertible linear map $\Phi$ 
satisfying 
$\Phi^{-1} = (d/|d|) \overline \Phi$ and 
$\trace(\Phi^*\Phi) = |d|$, we can recover
a unitary fiber functor. 

To rephrase the monoidal equivalence of fibre functors, 
consider an isomorphism of vector spaces $T: V \to W$. 
Given a rigidity pair $F(\epsilon): V\otimes V \to \C$ and 
$F(\delta): \C \to V\otimes V$, the composite maps 
\[
F(\epsilon)(T^{-1}\otimes T^{-1}): W \otimes W \to \C, 
\quad 
(T\otimes T)F(\delta): \C \to W\otimes W
\]
satisfy the rigidity relation. Hence, given another rigidity pair 
$G(\epsilon): W\otimes W \to \C$ and 
$G(\delta): \C \to W\otimes W$, 
$G(\epsilon) = F(\epsilon)(T^{-1}\otimes T^{-1})$ 
if and only if 
$G(\delta) = (T\otimes T)F(\delta)$. 
In terms of the associated linear maps 
$\Phi: V \to \overline V$ and $\Psi: W \to \overline W$, 
the condition is further equivalent to requiring 
$\Phi = {}^tT \Psi T$. 

Since 
\[
\Psi^{-1} = T\Phi^{-1}\,{}^tT, 
\qquad 
\overline \Phi = (T^*)^{-1}{\overline \Phi}\, {\overline T}^{-1}, 
% = \frac{d}{|d|} T\overline F {}^tT 
\]
the unitarity of $\Psi$ follows from 
that of $\Phi$ if $T$ is a unitary. 

To conclude the discussions so far, we introduce the following 
notation: Given a finite-dimensional Hilbert space $V$ and 
$d \in \R^\times = \R \setminus \{ 0\}$, set 
\[
\Gamma_d(V) = 
\left\{ 
\Phi: V \to \overline V; 
\Phi^{-1} = \frac{d}{|d|} \overline \Phi
\right\}.
\]

\begin{Lemma}
Let $\Phi \in \Gamma_d(V)$ and 
$\Psi \in \Gamma_d(W)$ with $V$ and $W$ 
finite-dimensional Hilbert spaces. 
Then these 
give rise to unitarily equivalent unitary fiber functors 
if and only if 
we can find a unitary map $T: V \to W$ 
satisfying $\Phi = {}^tT\Psi T$.

Conversely,  
given $\Phi \in \Gamma_d(V)$ and a unitary map $T: V \to W$, 
the linear map 
$\Psi = {}^tT^{-1} \Phi T^{-1}: W \to \overline W$ belongs to 
the set $\Gamma_d(W)$.
\end{Lemma}

Thus,
for each dimensionality of $V$ with $\dim V = n$, 
unitary fiber functors are classified up to unitary equivalence 
by the orbit space 
\[
\left. 
\Gamma_d(V)
\right/ U(V),
\]
where the unitary group $U(V)$ of $V$ acts on the set 
$\Gamma_d(V)$ by 
$\Phi.T = {}^tT \Phi T$ (a right action). 

We shall now rephrase the above orbit space 
by spectral data of $\Phi$. 
Let $\Phi \in \Gamma_d(V)$ and 
$\Phi = U|\Phi|$ be its polar decomposition 
with $U: V \to \overline V$ unitary and $|\Phi|: V \to V$ 
positive. 
By taking inverses
and substituting $\Phi^{-1}$ with $\pm \overline \Phi$, we have 
\[
U^{-1}U|\Phi|^{-1}U^{-1} = \pm \overline U\, \overline{|\Phi|}
\] 
and then, by the uniqueness of polar decompositions, 
\[
U^{-1} = \pm \overline U, 
\qquad 
U|\Phi|^{-1} U^{-1} = \overline{|\Phi|}.
\]

To utilize these identities, 
we introduce an antiunitary operator $C: V \to V$ by 
\[
Cv = \overline{Uv} = \overline U \overline{v}, 
\quad v \in V,
\]
which satisfies 
\[
C^2v = \overline{U\overline{U} \overline{v}} 
= \frac{d}{|d|} v, 
\]
i.e., $C^2 = \pm 1_V$. 
The equality 
$U|\Phi|^{-1} U^{-1} = \overline{|\Phi|}$ is then equivalent to 
the commutation relation 
\[
|\Phi|C = C |\Phi|^{-1}
\]
because of
\[
|\Phi|Cv = |\Phi|\overline{Uv} = 
\overline{\overline{|\Phi|} Uv}
= \overline{U|\Phi|^{-1}v} = C|\Phi|^{-1}v.
\]

Thus, the spectral structure of the positive operator $|\Phi|$ 
bears the symmetry of taking inverses on values 
under the operation of 
$C$; if $|\Phi|v = hv$ with $h>0$, then 
\[
|\Phi|Cv = C|\Phi|^{-1}v = h^{-1} Cv.
\]
In other words, 
$CV_h = V_{h^{-1}}$ if we denote the spectral subspace 
by $V_h = \{ v \in V; |\Phi|v = hv\}$. 
In particular, the antiunitary $C$ leaves 
$V_1$ invariant and here a dichotomy is in order
according to $C^2 = \pm 1_V$. 

For $C^2 = 1_V$, i.e., $d>0$, $C$ gives a real structure 
on $V$ and on the invariant subspace $V_1$ by restriction. 
Thus we can find an orthonormal basis $\{ w_k\}$ of $V_1$ 
which is real in the sense that $Cw_k = w_k$. 

For the case $C^2 = -1_V$ ($d<0$), 
the anti-unitarity of $C$ is used to see 
\[
(v|Cv) = (C^2v|Cv) = -(v|Cv), 
\]
i.e., $v \perp Cv$ for any $v \in V$. 
Hence we can find an orthonormal system $\{ w_k\}$ in $V_1$ such 
that $\{ w_k, Cw_k\}$ constitutes an orthonormal basis of $V_1$. 
Note that this occurs only when $V_1$ is even-dimensional. 

In both cases, the multiplicity information of the spectrum of 
$|\Phi|$ completely determines 
the operator $\Phi$ up to unitary equivalence:

\begin{Theorem}
The orbit space $\Gamma_d(V)/U(V)$ is completely parametrized 
by the eigenvalue list $\{ h_j\}$ of $|\Phi|$ 
(including the multiplicity but irrelevant to the order), 
satisfying $\{ h_j^{-1}\} = \{ h_j \}$, 
\[
\sum_j h_j^2 = |d|
\]
and $(d/|d|)^m = 1$, 
where $m$ denotes the multiplicity of an eigenvalue $1$
($m = \dim \ker(|\Phi| - 1_V)$). 
\end{Theorem}

\begin{Corollary}
The set $\Gamma_d(V)$ is empty unless $(d/|d|)^{\dim V} = 1$. 
\end{Corollary}

\begin{Example}
For a two-dimensional Hilbert space $V$ with an orthonormal basis 
$\{ v_1, v_2 \}$, the orbit space 
$\Gamma_d(V)/U(V)$ consists of one-point for each $\pm d \geq 2$. 
More concretely, 
with the expression $d = \pm(h^2 + h^{-2})$ ($h \geq 1$), 
the orbit is represented by a linear map $\Phi: V \to \overline V$, 
where
\[
\Phi(v_1) = h\overline{v_2}, 
\qquad 
\Phi(v_2) = \pm h^{-1} \overline{v_1}. 
\]

In other words, 
a unitary fiber functor $F: \cK_d \to \cH ilb$ such that 
$\dim F(X) = 2$ is unique up to unitary equivalence and realized by 
the choice
\begin{gather*}
F(\epsilon)(v_1\otimes v_2) = \pm h^{-1}, 
\quad 
F(\epsilon)(v_2\otimes v_1) = h,\\
F(\epsilon)(v_1\otimes v_1) = F(\epsilon)(v_2\otimes v_2) = 0.
\end{gather*}

With the notation in Example~2.3, the Hopf algebra is identified 
with $A_q$ for the choice $q = \mp h^{-2}$. 
The associated *-structure is then computed by the procedure in 
the appendix; let $\xi = (\xi_1,\xi_2)$ 
be the frame given by the orthonormal 
basis $\{ v_1, v_2\}$ 
and $\eta = (\eta_1,\eta_2)$ 
be the dual frame relative to $F(\epsilon)$. 
Then we have 
\[
(\xi_1,\xi_2) = (\eta_1,\eta_2)T 
\quad 
\text{with}\ 
T = 
\begin{pmatrix}
0 & h\\
\pm h^{-1} & 0
\end{pmatrix}
\]
and then 
\[
a^\eta = 
\begin{pmatrix}
a^* & b^*\\
c^* & d^*
\end{pmatrix}
= 
Ta^\xi T^{-1} = 
\begin{pmatrix}
0 & h\\
\pm h^{-1} & 0
\end{pmatrix}
\begin{pmatrix}
a & b\\
c & d
\end{pmatrix}
\begin{pmatrix}
0 & h\\
\pm h^{-1} & 0
\end{pmatrix}^{-1}, 
\]
i.e., 
\[
a^* = d,\qquad b^* = - q^{-1} c,
\]
which coincides with the choice in \cite{W1}. 
\end{Example}

% Thus the unitarity of the fiber functor $G$ is equivalent to 
% Since two unitary fiber functors $F$ and $G$ 
% on $\cK_d$ are unitarily equivalent if and only if 
% we can find a unitary map $T: V \to W$ which makes the following diagram 
% commutative
% \[
% \begin{CD}
% V\otimes V @>{T\otimes T}>> W\otimes W\\
% @V{F(\epsilon)}VV @VV{G(\epsilon)}V\\
% \C @= \C
% \end{CD}\quad ,
% \]
% the associated matrices $F$ and $G$ are related by $F = {}^tTGT$. 

\appendix
\section{}
We shall here review some points in
Tannaka-Krein duality on Hopf algebras following 
the basic references \cite{S, U} 
(cf.~\cite[\S 5.1]{CP} also for a concise description 
and \cite{Maj} for a further generalization)
%\cite{U} 
and present the reconstruction part in an informal way, 
which will be helpful in understanding the structure of 
relevant Hopf algebras 
%for non-experts in category theory like the author. 
without going into deep formalities of category theory.

\subsection{Saavedra Rivano-Ulblich's Theorem}

Let $\cC$ be a $\C$-linear category and 
consider a faithful linear functor $F: \cC \to \cV ec$, 
where $\cV ec$ is the tensor 
category of finite-dimensional $\C$-vector spaces. 
By imbedding $\cV ec$ into the tensor category $\overline{\cV ec}$ of 
(not necessarily finite-dimensional) $\C$-vector spaces, 
we regard $F$ as a functor from $\cC$ into $\overline{\cV ec}$. 

If we assume that the functor $F$ is exact, then we can define 
a $\C$-coalgebra $A$
% By the semisimplicity assumption on $\cC$, $F$ is automatically exact 
% and we obtain a $\C$-coalgebra $A$ 
as a solution of the universality problem for
natural transformations $F \to F\otimes B$, where $F\otimes B$ for an object 
$B$ in $\overline{\cV ec}$ is the functor specified by 
$(F\otimes B)(X) = F(X)\otimes B$ and $(F\otimes B)(f) = F(f)\otimes 1_B$ 
for $f \in \Hom(X,Y)$; a vector space $A$ together with a natural 
transformation $\rho: F \to F\otimes A$ is a solution of universality
in the sense that,
given a vector space $B$ and a natural transformation 
$\sigma: F \to F\otimes B$, we can find a unique linear map 
$\varphi: A \to B$ making the following diagram commutative 
for any object $X$. 
\[
\begin{CD}
F(X) @>{\rho_X}>> F(X)\otimes A\\
@| @VV{1_{F(X)}\otimes \varphi}V\\
F(X) @>>{\sigma_X}> F(X)\otimes B
\end{CD}\qquad .
\]

The coproduct $\Delta: A \to A\otimes A$ and 
the counit $\epsilon: A \to \C$ are then specified by the commutativity 
of the following diagrams for each $X$:
\[
\begin{CD}
F(X) @>{\rho_X}>> F(X)\otimes A\\
@V{\rho_X}VV @VV{1_{F(X)}\otimes \Delta}V\\
F(X)\otimes A @>>{\rho_X\otimes 1_A}> F(X)\otimes A\otimes A
\end{CD}\qquad ,
\quad
\begin{CD}
F(X) @>{\rho_X}>> F(X)\otimes A\\
@| @VV{1_{F(X)}\otimes \epsilon}V\\
F(X) @= F(X)\otimes \C
\end{CD}\qquad .
\]

If $\cC$ is a tensor category with unit object $A$ 
and $F$ is further assumed to be monoidal, 
then $A$ turns out to be a bialgebra with the product $m: A\otimes A \to A$ 
specified by the commutativity of the diagram 
\[
\begin{CD}
F(X\otimes Y) @>{\rho_{X\otimes Y}}>> F(X\otimes Y)\otimes A 
@>>> F(X)\otimes F(Y)\otimes A\\
@VVV @. @AA{1\otimes m}A\\
F(X)\otimes F(Y) @>>{\rho_X\otimes \rho_Y}> 
F(X)\otimes A\otimes F(Y)\otimes A @>>\text{flip}> 
F(X)\otimes F(Y)\otimes A\otimes A 
\end{CD}
\]
and the unit $1_A \in A$ given by the relation 
$\rho_I(1) = 1\otimes 1_A$ 
(recall that $F(I) = \C$ and $\rho_I: \C \to \C\otimes A$). 

The Ulbrich's theorem states that the constructed bialgebra $A$ is 
a Hopf algebra if and only if  
the tensor category $\cC$ is rigid.

\subsection{From Tensor Categories to Hopf Algebras}

We shall now present a naive (and less formal) approach 
to the results described above.
Start with an essentially small linear category $\cC$ and 
a faithful linear functor 
$F: \cC \to \cV ec$ again. 
Given an object $X$ of $\cC$, we denote a linear basis for 
the vector space $F(X)$ by the corresponding greek letter such as $\xi$, 
which is referred to as a frame. We also use the notation $|\xi|$ to 
stand for the object $X$. 
Let $\xi = (\xi_1,\dots, \xi_m)$ 
and $\eta = (\eta_1,\dots, \eta_n)$ be frames 
with $X = |\xi|$ and $Y = |\eta|$. 
For a linear map $T: F(X) \to F(Y)$, the matrix $(t_{ji})$ 
is denoted by $[\eta|T|\xi]$, where 
\begin{equation}
T(\xi_i) = \sum_{j=1}^n t_{ji} \eta_j, 
\quad \text{i.e.,}\quad
T\xi = \eta\,[\eta|T|\xi].
\end{equation}
 
Given a frame $\xi = (\xi_1,\dots,\xi_m)$, 
we introduce symbols $a^\xi_{ij}$ and let $A$ be the vector space 
spanned by all these symbols with the relations 
(called the covariance condition)
\begin{equation}
\sum_{l=1}^n t_{li}a^\eta_{kl} = \sum_{j=1}^m t_{kj}a^\xi_{ji}\ , 
\qquad 
1\leq i \leq m,\ 1 \leq k \leq n
\end{equation}
for various $X$, $Y$ and $T \in \Hom(X,Y)$. 
Note here that a basis-change $\sum_j t_{ji} \xi'_j = \xi_i$ 
induces the relation 
\begin{equation}
\sum_j a^{\xi'}_{ij} t_{jk} = 
\sum_j t_{ij} a^\xi_{jk}
\end{equation}
since $[\xi'|1_X|\xi] = (t_{ij})$.

More precisely, we consider a free vector space generated by symbols 
${\mathbf a}^\xi_{ij}$ and $A$ is defined to be the quotient space 
by the linear subspace spanned by 
(we have restricted ourselves to essentially small categories to 
avoid set-theoretical difficulties)
\begin{equation}
\sum_{l=1}^n t_{li}{\mathbf a}^\eta_{kl} - 
\sum_{j=1}^m t_{kj} {\mathbf a}^\xi_{ji}\,. 
\end{equation}
The symbol $a^\xi_{ij}$ then is used to denote the quotient element
of ${\mathbf a}^\xi_{ij}$. 

If we denote by $a^\xi$ 
the matrix arrangement of elements $a^\xi_{ij}$ 
in $A$ ($a^\xi \in M_m(A)$), 
then the covariance condition takes the form 
of matrix equations 
\begin{equation}
[\eta|T|\xi] a^\xi = a^\eta [\eta|T|\xi], 
\qquad 
T: X \to Y.
\end{equation}

A natural transformation $\rho: F \to F\otimes A$ is now defined so that 
\begin{equation}
\rho_X(\xi_i) = \sum_{j=1}^m \xi_j \otimes a^\xi_{ji}, 
\quad 1 \leq i \leq m, 
\ X = |\xi|,\ \xi = (\xi_1,\dots, \xi_m).
\end{equation}
The well-definedness of $\rho$ is exactly the covariance condition.  
The natural transformation $\rho$ is a solution to 
the universality problem. In fact, given another natural transformation 
$\sigma: F \to F\otimes B$, 
the linear map $\varphi: A \to B$ is well-defined 
by $\varphi(a^\xi_{ji}) = b_{ji}$ with
\begin{equation}
\sigma_X(\xi_i) = \sum_{j=1}^m \xi_j\otimes b_{ji}.
\end{equation}

The coproduct $\Delta: A \to A\otimes A$ is well-defined by 
\begin{equation}
\Delta(a^\xi_{ij}) = \sum_{k=1}^m a^\xi_{ik}\otimes a^\xi_{kj}
\end{equation}
because it preserves the covariance condition; 
\begin{align*}
\sum_{l=1}^n t_{li} \sum_{s=1}^n &a^\eta_{ks}\otimes a^\eta_{sl} 
- \sum_{j=1}^m t_{kj} \sum_{r=1}^m a^\xi_{jr}\otimes a^\xi_{ri}\\
&= 
\sum_{s=1}^n a^\eta_{ks}\otimes \sum_{r=1}^m t_{sr}a^\xi_{ri} 
- \sum_{r=1}^m \sum_{s=1}^n t_{sr}a^\eta_{ks}\otimes a^\xi_{ri}
= 0.
\end{align*} 

The coproduct $\Delta$ is obviously coassociative by its form 
and the counit $\epsilon: A \to \C$ for $\Delta$ is also well-defined by 
$\epsilon(a^\xi_{ij}) = \delta_{ij}$
simply because of 
\begin{equation}
\sum_{l=1}^n t_{li}\delta_{kl} - \sum_{j=1}^m t_{kj}\delta_{ji} 
= 0.
\end{equation}

The coproduct is chosen so that each $\rho_X: F(X) \to F(X)\otimes A$ 
gives a corepresntation of $A$. The resulting right comodule is denoted by 
$F(X)^A$. In this way, 
we have obtained a coalgebra $A$ and the functor $F$ is lifted to 
a faithful linear functor
of $\cC$ into the category $\cM^A$ 
of finite-dimensional right $A$-comodules. 
(Note that the $A$-colinearity is exactly the covariance 
condition.)

Non-trivial is the fact that the functor $\cC \to \cM^A$ gives 
an equivalence of categories if and only if $\cC$ is abelian and 
$F$ is exact (\cite[Theorem 2.3.5]{S}). 
For a semisimple $\cC$, however, the equivalence theorem is 
an easy consequence of irreducible decompositions of objects 
as it will be discussed below.

Now assume that $\cC$ is a tensor category with unit object $I$ and 
$F$ is monoidal. We may assume $F$ is strictly monoidal 
without loss of generality. 
Then $A$ is an algebra with the multiplication map 
$\mu: A\otimes A \to A$ defined by 
\begin{equation}
a^\xi_{ij}a^\eta_{kl} = a^{\xi\otimes\eta}_{i,k;j,l},
\end{equation}
where $\xi\otimes\eta = \{ \xi_i\otimes \eta_k \}$ is a frame of $X\otimes Y$ 
obtained from $\xi$ and $\eta$ by taking tensor products. 
This is again well-defined because it preserves the covariance condition; 
the relation 
\begin{equation}
\left(
\sum_l a^\eta_{kl}t_{li} - \sum_j t_{kj}a^\xi_{ji}
\right)
a^\zeta_{rs}
= \sum_l a^{\eta\otimes\zeta}_{k,r;l,s} t_{li} 
- \sum_j t_{kj} a^{\xi\otimes\zeta}_{j,r;i,s} 
\end{equation}
is, for example, associated to the morphism 
$T\otimes 1_Z: X\otimes Z \to Y\otimes Z$ 
in $\cC$.

Clearly $\mu$ is associative and 
the coproduct $\Delta$ is compatible 
with $\mu$; 
\begin{equation}
\Delta(a^{\xi\otimes\eta}_{i,k;j,l}) 
= \sum_{r,s} a^{\xi\otimes\eta}_{i,k;r,s}\otimes 
a^{\xi\otimes \eta}_{r,s;j,l}
= \sum_{r,s} a^\xi_{ir}a^\eta_{ks}\otimes a^\xi_{rj}a^\eta_{sl}
= \Delta(a^\xi_{ij})\Delta(a^\eta_{kl}).
\end{equation}
The unit $1_A$ for $\mu$ is given by $\rho_I(1) \in \C\otimes A = A$ 
($F(I) = \C$). 

At this point, we have established a bialgebra structure on $A$ 
so that $\cC \to \cM^A$ is a monoidal imbedding. 

Now assume the rigidity on $\cC$ and we shall show that 
$A$ is a Hopf algebra. 
Let $\xi$ be a frame for $X$ and choose a dual object $X^*$ together 
with a non-degenerate morphism $\epsilon_X: X^*\otimes X \to I$ 
($\epsilon_X$ is said to be non-degenerate if we can find 
a morphism $\delta_X: I \to X\otimes X^*$ fulfilling  
$(1_X\otimes \epsilon_X)(\delta_X\otimes 1_X) = 1_X$ and 
$(\epsilon_X\otimes 1_{X^*})(1_{X^*}\otimes \delta_X) = 1_{X^*}$). 
We then define a frame $\xi^*$ for $X^*$ by the relation
\begin{equation}
F(\epsilon_X)(\xi^*_i\otimes \xi_j) = \delta_{i,j},
\end{equation}
i.e., $\xi^*$ is the dual basis of $\xi$ with respect to 
the non-degenerate bilinear form 
$F(\epsilon_X): F(X^*)\otimes F(X) \to \C$. 
By the covariance condition, 
we see that the matricial element $a^{\xi^*} \in M_m(A)$ does not 
depend on the choice of either $X^*$ or $\epsilon_X$. 
The linear isomorphism $S: A \to A$ is then well-defined by 
$S(a^\xi_{ij}) = a^{\xi^*}_{ji}$. The property of antipode 
\begin{equation}
\sum_r S(a^\xi_{ir}) a^\xi_{rj} = 
\epsilon(a^\xi_{ij}) 1_A 
= \sum_r a^\xi_{ir} S(a^\xi_{rj})
\end{equation}
follows from the commutativity of diagram
\begin{equation}
\begin{CD}
F(X^*)\otimes F(X) @>{\rho_{X^*\otimes X}}>> 
F(X^*)\otimes F(X)\otimes A\\
@V{F(\epsilon_X)}VV @VV{F(\epsilon_X)\otimes 1_A}V\\
\C @>>{\rho_I}> A\\
@V{F(\delta_X)}VV @VV{F(\delta_X)\otimes 1_A}V\\
F(X)\otimes F(X^*) @>>{\rho_{X\otimes X^*}}> 
F(X)\otimes F(X^*)\otimes A
\end{CD}\quad .
\end{equation}
Notice here that the linear map $F(\delta_X): \C \to F(X)\otimes F(X^*)$ 
takes the form 
\begin{equation}
F(\delta_X): 1 \mapsto \sum_{j=1}^m \xi_j\otimes \xi^*_j
\end{equation}
by the uniqueness of the coparing map in $\cV ec$. 

As a conclusion of discussions so far, we have 

\begin{Proposition}
Let $\cC$ be an essentially small rigid tensor category and 
$F: \cC \to \cV ec$ be a faithful monoidal functor from 
$\cC$ to the tensor category $\cV ec$ of 
finite-dimensional vector spaces. 
Then we can find a Hopf algebra $A$ and a faithful monoidal functor 
from $\cC$ to the tensor category $\cM^A$ of finite-dimensional right 
$A$-comodules so that $F$ is the forgetful functor of this.
\end{Proposition}

We shall now assume that the linear category $\cC$ is 
semisimple in the sense that  
any object is isomorphic to 
a direct sum of simple objects, where an object $X$ in $\cC$ 
is said to be simple if $\End(X) = \Hom(X,X) = \C 1_X$. 
Let $S$ be the set of isomorphism classes of 
simple objetcs of $\cC$ and choose a representative family 
$\{ X_s \}_{s \in S}$ of simple objects. 
Furthermore, select a frame $\xi$ for each $X_s$ and denote 
the associated element of $A$ by $a^{(s)}_{ij}$.  
% It would be worth pointing that the triviality of end-algebra is
% not sufficient to imply the simplicity if we deal with non-semisimple 
% categories. 
% Anyway, all categories are asuumed to satisfy the above semisimplicity 
% in the rest of this appendix. 

Then by decomposing objects into direct sums of $X_s$'s, 
we see that the family $\{ a^{(s)}_{ij} \}_{s, i, j}$ 
constitutes a linear basis of $A$. 
Observing the coproduct formula 
\begin{equation}
\Delta(a^{(s)}_{ij}) = \sum_k a^{(s)}_{ik}\otimes a^{(s)}_{kj},
\end{equation}
this implies that the coalgebra $A$ is semisimple 
with all simple comodules supplied by $F(X_s)^A$ ($s \in S$).
Thus we have arrived at 

\begin{Corollary}
Assume that $\cC$ is semisimple furthermore. 
Then $\cC$ is monoidally equivalent to 
the tensor category $\cM^A$.
\end{Corollary}

Assume that the (strict) tensor category $\cC$ is 
generated by an object $X$; objects of $\cC$ are of the form 
$X^{\otimes n} = X\otimes \dots \otimes X$ (the $n$-th tensor 
power of $X$) for $n=0, 1, 2, \dots$. Here $X^{\otimes 0}$ is 
the unit object $I$ by definition. 
Choose a frame $\xi$ of $X$ once for all and write 
$a_{ij} = a^\xi_{ij}$. 

Let $F: \cC \to \cV ec$ be a faithful strictly monoidal functor and 
$A$ be the associated bialgebra. 
Then $A$ is generated by the elements $a_{ij}$ as a unital algebra 
with the relations given by the commutativity of diagrams 
\begin{equation}
\begin{CD}
F(X^m) @>{\rho_m}>> F(X^m)\otimes A\\
@V{F(f)}VV @VV{F(f)\otimes 1_A}V\\
F(X^n) @>>{\rho_n}> F(X^n)\otimes A 
\end{CD}
\end{equation}
for various morphisms $f: X^m \to X^n$, 
with the coproduct $\Delta: A \to A\otimes A$ specified by 
\begin{equation}
\Delta(a_{ij}) = \sum_k a_{ik}\otimes a_{kj}
\end{equation}
as already discussed. 

Now we restrict ourselves to the Temperley-Lieb category 
$\cK_d$. In this case, $\Hom(X^m,X^n)$ is generated by 
basic arcs as tensor cateogries, which implies that the 
algebrta $A$ is generated by $a_{ij}$ with the relations 
given by the commutativity of two diagrams 
\begin{equation}
\begin{CD}
F(X)\otimes F(X) @>{\rho_{X\otimes X}}>> 
F(X)\otimes F(X)\otimes A\\
@V{F(\epsilon_X)}VV @VV{F(\epsilon_X)\otimes 1_A}V\\
\C @>>> \C\otimes A\\
@V{F(\delta_X)}VV @VV{F(\delta_X)\otimes 1_A}V\\
F(X)\otimes F(X) @>>{\rho_{X\otimes X}}> 
F(X)\otimes F(X)\otimes A
\end{CD}\ .
\end{equation}

If we write 
\begin{equation}
E_{ij} = F(\epsilon_X)(\xi_i\otimes \xi_j), 
\qquad 
F(\delta_X)(1) = \sum_{ij} D_{ij} \xi_i\otimes \xi_j, 
\end{equation}
then $D = (D_{ij})$ is the inverse matrix of $E = (E_{ij})$ 
and the above commutativity takes the form
\begin{equation}
\sum_{k,l} E_{k,l} a_{ki} a_{lj} = E_{ij} 1_A,
\qquad 
\sum_{i,j} D_{ij} a_{ki} a_{lj} = D_{kl} 1_A
\end{equation}
of quadratic relations. 

\subsection{From Unitary Fiber Functors to Compact Quantum Groups}

From here on, $\cC$ is assumed to be a rigid C*-tensor category and 
$F: \cC \to \cH ilb$ be a unitary fiber functor 
(see \cite{W2, URT, MRT}). 
The associated Hopf algebra $A$ is then a Hopf *-algebra. 
In fact, let $\xi$ be a frame consisting of orthonormal vectors and 
define a conjugate-linear map $A \ni a \mapsto a^* \in A$ so that 
$(a^\xi_{ij})^* = S(a^{\xi}_{ji}) = a^{\xi^*}_{ij}$. 
This is well-defined because it preserves the covariance condition 
as seen from 
\begin{equation}
\sum_{l=1}^n \overline{t_{li}}(a^\eta_{kl})^* 
- \sum_{j=1}^m \overline{t_{kj}}(a^\xi_{ji})^*
= S\left( 
\sum_{l=1}^n \overline{t_{li}} a^\eta_{lk} 
- \sum_{j=1}^m \overline{t_{kj}} a^\xi_{ij}
\right), 
\end{equation}
where 
\begin{equation}
\sum_{l=1}^n \overline{t_{li}} a^\eta_{lk} 
- \sum_{j=1}^m \overline{t_{kj}} a^\xi_{ij}
\end{equation}
vanishes as the covariance condition for 
the adjoint morphism $T^*: Y \to X$.

We next check that the operation * is involutive. 
To this end, choose an orthonormal frame $\eta$ for $X^*$ and 
write 
\begin{equation}
\xi^*_i = \sum_j (\eta_j|\xi^*_i) \eta_j.
\end{equation}
This change-of-basis gives rise to the relation 
\begin{equation}
\sum_k a^\eta_{ik} (\eta_k|\xi^*_j) 
= \sum_k (\eta_i|\xi^*_k) a^{\xi^*}_{kj}
\end{equation}
and then, by applying the *-operation, 
\begin{equation}
\sum_k (\xi^*_j|\eta_k)\, a^{\eta^*}_{ik} 
= \sum_k (\xi^*_k|\eta_i) \left( a^{\xi^*}_{kj} \right)^*.
\label{star}
\end{equation}

For an explicit computation of $\eta^*$, we use 
$F(\delta_X)^*: F(X^*)\otimes F(X) \to\C$ 
as a non-degeneate pairing; 
$F(\delta_X)^*(\eta^*_k\otimes \eta_l) = \delta_{k,l}$. 
Since $F(\delta_X) = \sum_i \xi_i\otimes \xi^*_i$, this means 
\begin{equation}
\sum_j (\xi_j|\eta^*_l) (\xi^*_j|\eta_k) = \delta_{k,l}, 
\quad \text{i.e.,}\ 
\sum_k (\xi_i^*|\eta_k)(\xi_j|\eta^*_k) = \delta_{i,j}. 
\end{equation}
which is used in (\ref{star}) to obtain 
\begin{equation}
\left( 
a^{\xi^*}_{lj}
\right)^* 
= \sum_{i,k} (\xi_j^*|\eta_k) (\xi_l|\eta^*_i) 
a^{\eta^*}_{ik}.
\end{equation}
If we use the change-of-basis 
$\eta^*_i = \sum_l (\xi_l|\eta^*_i) \xi_l$
in the last expression, it takes the form 
\begin{equation}
\sum_{i,k} (\xi^*_j|\eta_k) (\xi_i|\eta^*_k) a^\xi_{li}
= \sum_i \delta_{i,j} a^\xi_{li} = a^\xi_{lj},
\end{equation}
showing $((a^\xi_{lj})^*)^* = (a^{\xi^*}_{lj})^* = a^\xi_{lj}$. 

\begin{Remark}
The dual frame $\xi^*$ is not necessarily orthonormal 
and we can not expect the relation such as $\xi^{**} = \xi$. 
\end{Remark}

It is now immediate to check that the involution $*$ makes 
$A$ a Hopf *-algebra. In fact, given orthonormal frames 
$\xi$ and $\eta$, their tensor product frame $\xi\otimes \eta$ is 
orthonormal and we have 
\begin{equation}
(a^\xi_{ij}a^\eta_{kl})^* 
= (a^{\xi\otimes \eta}_{i,k;j,l})^* 
= S(a^{\xi\otimes \eta}_{j,l;i,k}) 
= S(a^\xi_{ji}a^\eta_{lk})
= S(a^\eta_{lk}) S(a^\xi_{ji})
= (a^\eta_{kl})^* (a^\xi_{ij})^*,
\end{equation}
whereas the compatibility with the coproduct is seen from 
\begin{equation}
\Delta((a^\xi_{ij})^*) 
= \Delta(a^{\xi^*}_{ij}) 
= \sum_k a^{\xi^*}_{ik}\otimes a^{\xi^*}_{kj} 
= \sum_k (a^\xi_{ik})^*\otimes (a^\xi_{kj})^* 
= (\Delta(a^\xi_{ij}))^*.
\end{equation}

Finally observe that each $A$-comodule 
$\rho_X: F(X) \to F(X)\otimes A$ is unitary in the sense 
(see \cite{W1}, cf.~also \cite{MMNNU}) that 
\begin{equation}
\langle \rho_X(v)|\rho_X(w)\rangle 
= (v|w) 1_A
\quad
\text{for $v$, $w \in F(X)$},
\end{equation}
where $\langle\ |\ \rangle$ on the left side 
is the $A$-valued inner product defined by 
\begin{equation}
\langle v\otimes a|w\otimes b\rangle = (v|w) a^*b,
\quad a, b \in A.
\end{equation}
This follows from 
\begin{align*}
\langle \rho_X(\xi_i)|\rho_X(\xi_j)\rangle 
&= \left\langle 
\sum_k \xi_k\otimes a^\xi_{ki} \left| 
\sum_l \xi_l\otimes a^\xi_{lj} \right.
\right\rangle\\
&= \sum_{k,l} 
(\xi_k|\xi_l) (a^\xi_{ki})^* a^\xi_{lj}\\
&= \sum_k (a^\xi_{ki})^* a^\xi_{kj} = 
\sum_k S(a^\xi_{ik}) a^\xi_{kj} 
= \delta_{ij} 1_A.
\end{align*}

\begin{Remark}
In the definition of compatible *-operations, 
we do not use the positivity of inner products and 
can equally well work with *-tensor categories and 
(finite-dimensional) indefinite inner product spaces, 
which in fact arizes from the Temperley-Lieb category 
$\cK_d$ for $-2 < d <2$ 
(the associated Hopf *-algebra is then $SU_q(1,1)$ 
with $|q| = 1$, another real form of $SL_q(2,\C)$). 
\end{Remark}

Thus we have derived the algebraic part in 
Woronowicz' Tannaka-Krein duality (\cite{W2}). 

\begin{Proposition}
Let $\cC$ be a rigid abelian C*-tensor category with simple unit 
object and $F: \cC \to \cH ilb$ be a unitary fiber functor. 
Then we can find a *-Hopf algebra so that $\cC$ is C*-monoidally 
equivalent to the C*-tensor category of finite-dimensional 
unitary $A$-comodules with $F$ identified with the associated 
forgetful functor. 
\end{Proposition}

\end{document}